\documentclass[a4,12pt,multicol,makeidx]{article}
\def\origin{
  \clearpage
\vskip-\baselineskip\vskip-\topskip%
  \vbox to 0pt{\vskip-1in%
    \hbox to 0pt{\hskip-1in%
      \hbox to 0pt{\vrule width 1cm height .4pt depth 0mm\hss}%
      \vbox to 0pt{\hrule width .4pt height 0pt depth 1cm\vss}%
    \hss}%
  \vss}
  \vskip-\baselineskip
  \vbox to 0pt{\vskip-1in\vskip3cm%
    \hbox to 0pt{\hskip-1in\hskip3cm%
      \hbox to 0pt{\hss\vrule width 2cm height .4pt depth 0mm\hss}%
      \vbox to 0pt{\vss\hrule width .4pt height 1cm depth 1cm\vss}%
    \hss}%
  \vss}%
\vskip5mm\hskip10mm (3cm,3cm)
}%

 \def\l{\lambda}   \def\p{\partial}  \def\e{\varepsilon} 
    
\def\leq{\underline{<}}

\newenvironment{theorem}{%
\par \bigskip \it}{%
\bigskip \par}

\makeindex
\title{Extended groups of semigroups and  backward problems of heat equations
}
\author{Mariko Arisawa\\ Department of Mathematics,
\\Ochanomizu University\\
2-1-3 Otsuka Bunkyo-ku, Tokyo, Japan\\
}
\date{}
\begin{document}
\maketitle
\bigskip

\section{Introduction - the origin of the problem.} 

We consider the heat equation in $[0,1]\subset{\bf R}$ :
\begin{eqnarray}
	\frac{\p u}{\p t}(t,x)&=&\Delta u(t,x) 
	\qquad t\geq 0,\quad x\in [0,1], \nonumber\\
	u(0,x)&=&u_0(x) \quad x\in [0,1],\\
	u(t,0)&=&u(t,1)=0 \quad t\geq 0,\nonumber
\end{eqnarray}
where the initial condition 
$u_0(x)$ belongs to $L^2(0,1)=H$. Let $\{{\bf e_n }\}$ be a complete orthonormal basis of $L^2(0,1)$. Then, for the initial condition $u_0(x)$$=\sum_{n=1}^{\infty} a_n {\bf e_n}$, the solution $u(t,x)$ is given by the semigroup $\{T_t\}$ :
$$
	T_t u_0(x)= \sum_{n=1}^{\infty} a_n \exp(-n^2 \pi^2 t){\bf e_n}, 
$$
which traces out the  trajectory of the solution of (1) on $H$. \\

If we go along the trajectories of solutions to the direction negative in time, we remark that there are two kinds of situations : (i) the trajectory which continues until $t=-\infty$; and (ii) the trajectory which  continues backward  only for finite period of time. The initial condition $u_0(x)$ of the latter trajectory, which stops at finite $t\leq 0$ can be written 
$u_0(x)=T_t u_1(x)$, where $u_1(x)$ is the initial condition  from which the backward trajectory does not exist. Thus, the case of (ii) can be studied by the initial value  $u_1$, at which the solution trajectory is irreversible to the past. 
For the case of (ii), we can characterize the initial function $u_0(x)$$=\sum_{n=1}^{\infty} a_n{\bf e_n}$ from which the solution trajectory goes back to $t=-\infty$, by : 
$$
	D=\{ \sum_{n=1}^{\infty} a_n{\bf e_n}| \quad
		\sum_{n=1}^{\infty} a_n^2 \exp (2n^2{\pi}^2 t) <\infty\quad 
		\hbox{for}\quad \forall t\geq 0
	\}.
$$
On the other hand, the initial function $u_0(x)$$=\sum_{n=1}^{\infty} a_n{\bf e_n}$ from which the solution trajectory cannot be defined 
 for any $t<0$  satisfies : 
$$
	Z=\{ \sum_{n=1}^{\infty} a_n{\bf e_n}| \quad
		\sum_{n=1}^{\infty} a_n^2 <\infty,\quad 
		\sum_{n=1}^{\infty} a_n^2 \exp (2n^2{\pi}^2 t)=\infty\quad 
		\hbox{for}\quad \forall t\geq 0
	\}. 
$$
If $u_0\in D$, then
$$
	T_{-t} u_0(x)=\sum_{n=1}^{\infty} a_n \exp (n^2{\pi}^2 t){\bf e_n}
	\in H=L^2(0,1).
$$
For $u_0\in Z$, though we can formally write 
$$
	T_{-t} u_0(x)=\sum_{n=1}^{\infty} a_n \exp (n^2{\pi}^2 t){\bf e_n}, 
$$
 $T_{-t} u_0$ is not in $H$, and should belong to a functional space larger than $H=L^2(0,1)$. For example, for the following Hilbert space
$$
	H_t=\{\sum_{n=1}^{\infty} a_n{\bf e_n}|\quad
	\sum_{n=1}^{\infty} a_n^2 \exp (-2n^2{\pi}^2 t)<\infty
	\},
$$
equipped with the norm
$$
	\|\sum_{n=1}^{\infty} a_n \exp (n^2{\pi}^2 t)\|_t
	= \|\sum_{n=1}^{\infty} a_n \exp (-n^2{\pi}^2 t){\bf e_n}\|,
$$
where the $\|\cdot\|$ is the norm in $H$.
We have 
$$
	H\subset H_t,\quad T_{-t}u_0(x)\in H_t \quad \hbox{if}\quad u_0(x)\in Z. 
$$
Thus, the solution trajectory which stops at a point in the Hilbert space $H$ is extendable in the larger space $H_t$. Therefore, the heat semigroup $T_t$ defined by (1), is extendable backward in time, by considering the larger functional spaces equipped with the above kind of norms $\|\cdot\|_t$. This construction of the sequence of spaces $\{H_t\}$ comes from  the explicite expression :
$$
	T_t u_0(x)=\sum_{n=1}^{\infty} a_n \exp (\lambda_n t){\bf e_n}
	\quad \l_n<0,\quad t<0\quad
	(u_0(x)=\sum_{n=1}^{\infty} a_n {\bf e_n}),
$$
where $\lambda_n=-(n \pi)^2$, eigenvalues  of $-\Delta$ : 
$$
	-\Delta= \left(
	\begin{array}{ccccc}
	-\pi^2 & \quad & \quad & \quad & \quad\\
	\quad  & -(2\pi)^2& \quad & \quad & \quad\\
	\quad & \quad & \cdots & \quad & \quad\\
	\quad & \quad & \quad & -(n \pi)^2 & \quad\\
	\quad & \quad & \quad & \quad & \cdots
	\end{array}
	\right). 
$$
The above consideration leads to the following question. If an infinitesmal generator $A$ of a semigroup $T_t$ in a functional space $X$ does not have the spectrum decomposition as above, is there a larger functional space $X_t$ in which any solution trajectory which starts from any point in $X$ continues backward in time uniformly until $-t<0$ ? Moreover, is there an extended space $E$ of $X$, in which any solution trajectory starting from any point in $E$, can go back to the past until $t=-\infty$ ?
 In other words, we are interested in the following question. \\
 
{\bf Problem} Is there an extended functional space in which a semigroup becomes a group ? 
Find the condition for such an extension. \\

In this paper, under the following two conditions (I) and (II), we show the existence of the functional space $E$ described in the above problem. \\
(I) The set 
$$
	D=\{x\in X|\quad \forall t>0\quad \exists y=y(t)\in X\quad \hbox{s.t.} \quad 
	T_ty=x
	\}
$$ 
is dense in $X$.\\
(II) The semigroup $T_t$ has the backward uniqueness property, i.e. 
$$
	\hbox{if}\quad T_t x= T_t y\quad \hbox{for some}\quad t>0, \quad
	\hbox{then}\quad x=y.
$$
The following is the plan of this paper : in \S 2, a sufficient condition such that $D$ satisfies (1) is given, as well as some properties of $D$; in \S 3, we present  the extended space $E$ explicitely; and in \S 4,  the structure of $E$ and the relationship between $E$ and $D$ are investigated. \\

This paper was written in 1989 as my master's thesis, by receiving valuable advises from Professor Yukio Komura. \\

\section{The set of points from which trajectories continue negative in time until $t=-\infty$.}

From this section, we denote
$$
	D=\{x\in X|\quad \forall t>0\quad \exists y=y(t)\in X\quad \hbox{s.t.}\quad     T_ty=x\},
$$
for the set of points from which trajectories continue negative in time until $t=-\infty$, and 
$$
	D_t=\{x\in X|\quad \exists y\in X\quad \hbox{s.t.}\quad     T_ty=x\},
$$
for the set of points from which trajectories continue negative in time at least until $-t<0$. \\

The following is a sufficient condition such that $D$ is dense in $X$.\\

{\bf Lemma 1.}
\begin{theorem}
Let $T_t$ ($t\geq 0$) be a linear semigroup defined on $X$. Then, the necessary and sufficient condition such that $D$ is dense in $X$ is :
$$
	\forall t>0 \quad D_t=T_tX \quad \hbox{is dense in }\quad X.
$$
\end{theorem}
$Proof.$ The necessary condition is clear. We only show the part of the sufficient condition. It is enough to prove for any $x_0\in X$ and for any $\e_0>0$, 
$$
	U_{\e_0}(x_0)\cap D \neq \emptyset.
$$
Since $D_1$ is dense in $X$, there exists $x_1\in X$ such that $T_1 x_1\in U_{\e_0}(x_0)$. By repeating this, for any $\e_k>0$ there exists $x_{k+1}\in X$ such that $T_1x_{k+1}\in U_{\e_k}(x_k)$. In particular, we can take $\{\e_n\}_{n\in {\bf N}}$ such that the sequence $\{T_n x_n\}$ becomes a Cauchy sequence in $X$, since 
$$
	||T_nx_n - T_{n-1}x_{n-1}||\leq 
	M_{n-1}||T_1 x_n -x_{n-1}||\leq M_{n-1}\e_{n-1},
$$
where $M_t>0$ is a number such that $||T_t||\leq M_t$, more precisely $M_t=M e^{wt}$ ($M>0$ a constant). For such a choice of $\{x_n\}$, there exists a limit
$$
	\lim_{n\to \infty} T_n x_n= \exists x_{\infty,0}. 
$$
Moreover, for each $k$, since 
$$
	||T_{n-k}x_n - T_{n-1-k}x_{n-1}||\leq 
	M_{n-1-k}||T_1 x_n -x_{n-1}||\leq M e^{-kw}e^{(n-1)w}\e_n,
$$
 there exists a limit 
$$
	\lim_{n\to \infty} T_{n-k} x_n= \exists x_{\infty,k}. 
$$
Thus,
$$
	T_kx_{\infty,k}= \lim_{n\to \infty} T_{n} x_n=x_{\infty,0},
$$
and $x_{\infty,0}\in D$. 
Since
$$
	||x_{\infty,0}-x_0||=\lim_{n\to \infty}||T_nx_n-x_0||,
$$
and
$$
	||T_nx_n - x_0||\leq ||T_{n}x_n - T_{n-1}x_{n-1}||+ 
	||T_{n-1}x_{n-1} - T_{n-2}x_{n-2}||+\cdots + ||T_{1}x_1 - x_{0}||
	\leq \sum_n Me^{wt}\e_n,
$$
by choosing $\e_n>0$ appropriately, we have $||x_{\infty,0}-x_0||<\e_0$,
 which proves the claim. \\

 Lemma 1 holds also for the nonlinear semigroup $T_t$. Since, in general, the nonlinear semigroup is defined on $\overline{D(A)}$ ($A$ is the infinitesmal generator of $T_t$), not necessarily on whole $X$, we rewrite the lemma as follows.\\
 
{\bf Lemma 2.}
\begin{theorem}
Let $T_t$ be a semigroup (possibly nonlinear) defined on $\overline{D(A)}$, satisfying
$$
	\|T_tx-T_ty\|\leq M e^{wt}\|x-y\| \quad \forall x,y\in \overline{D(A)},
$$
where $M$, $w$ are positive constants independent on $t>0$. Then, the necessary and the sufficient condition that $D$ is dense in $\overline{D(A)}$ is 
 that $D_t$ is dense in $\overline{D(A)}$ for any $t>0$. 
\end{theorem}
 
$Proof.$ As before, it is enough to prove the part of the sufficient condition.  As in the proof of Lemma 1, for any $x_0\in \overline{D(A)}$ and for any $\e_0>0$, 
  we can take inductively a sequence of numbers $\e_k>0$ ($k\geq 1$), and a sequence of points $x_k\in \overline{D(A)}$ ($k\geq 1$) such that 
 $$
 	T_1x_{k+1}\in U_{\e_k}(x_k),
 	\quad \lim_{n\to \infty}T_nx_n=\exists x_{\infty,0}\in D,
 $$
 where $\|x_{\infty,0}-x_0\|\leq \e_0$. This proves the claim. \\
 
 By using Lemmas 1 and 2, we can now give sufficient conditions such that $D$, the set of points in $\overline{D(A)}$ from which the trajectories continue until $t=-\infty$, is dense in $\overline{D(A)}$. \\
 
{\bf Proposition 1.}
\begin{theorem}
	Let $T_t$ be a linear holomorphic semigroup on a Banach space $X$. That is, for any $t_0>0$, there exists a neighborhood $U$ of $t_0$, and the following  Taylor development :
$$
	T_t=\sum_{k=0}^{\infty} (t-t_0)^k A_k \quad t\in U.
$$
Then, the set 
$$
	D=\{x\in X|\quad \forall t>0 \quad \exists y=y(t)\quad \hbox{s.t.}\quad 
	T_ty=x\}
$$
is dense in $X$. 
\end{theorem}
$Proof.$ Put
$$
	D_{t_1}=\{x\in X|\quad \exists y\in X \quad T_{t_1}y=x\}.
$$
We use the argument by contradiction. 
 From Lemma 1, we assume that there exists $t_1>0$ such that $D_{t_1}$ is not dense in $X$, and we shall lead to a contradiction. If $\overline{D_{t_1}}$, closed and convex, is not $X$, by the Hahn-Banach theorem there exists $f\in X'$ (dual space of $X$), such that 
 $$
 	f\neq 0,\quad <x,f>=0 \quad \forall x\in \overline{D_{t_1}}.
$$ 
 Since $T_t$ is holomorphic, there exists a neighborhood $U$ od $t_1$ such that \begin{equation}\label{holo}
 	<T_tx,f>=<\sum_{k=0}^{\infty}(t-t_1)^k A_kx,f>
 	=\sum_{k=0}^{\infty}(t-t_1)^k <A_kx,f>.
\end{equation}

For $t>t_1$, $T_tx\in D_{t_1}$, and $<T_tx,f>=0$. 
 Since the right-hand side of (\ref{holo}) is a holomorphic function, from  the unique continuation theorem, 
 $$
 	<T_tx,f>=0 \quad \forall t\in U. 
 $$
 By repeating this argument, we have 
 $$
 	<T_tx,f>=0 \quad \forall t>0,\quad \forall x\in X.
 $$
 Moreover, from
 $$
 	<x,f>=\lim_{t\downarrow 0} <T_t x,f>=0 \quad \forall x\in X,
 $$
 we get $f=0$. This contradicts to our previous assumption, and thus we have proved the claim. \\
 
{\bf Corollary 1.}
\begin{theorem}
Let $A$ be an infinitesmal generator of a linear holomorphic semigroup. Consider the following evolution equation in $X$
$$
	\frac{dx}{dt}+Ax=f(t)\quad x(0)=x_0,
$$
where $f$ is a locally H{\"o}lder continuous function from $(0,\infty)$ to $X$, such that 
$$
	\int_0^{\rho} \|f(t)\|dt<+\infty\quad \hbox{for some}\quad \rho>0,
$$
and assume that there exists a global solution $x(t)$ :
$$
	x(t)=T_t x_0+ \int_0^t T_{t-s}f(s)ds.
$$
Let $S_t$ be the semigroup 
$$
	S_t x_0=x(t).
$$ 
Then, the set 
$$
	D=\{x\in X|\quad \forall t>0 \quad \exists y=y(t)\quad \hbox{s.t.}\quad S_ty=x\}
$$
is dense in $X$.
\end{theorem}
$Proof.$ From Lemma 2, it is enough to show that for any $t>0$, the set 
$$
	D_t=\{x\in X|\quad \exists y\in X \quad \hbox{s.t.}\quad 
	T_t y+ \int_0^t T_{t-s}f(s)ds=x \}
$$
is dense in $X$. However, it is clear, because 
$$
	D_t=T_t X + \int_0^t T_{t-s}f(s)ds, 
$$
 where $+$ is the parallel translation. \\
 
 \medskip
 
{\bf Remark.} The example of the heat equation in \$.1 is a special case of Corollary 1. \\
 
 \section{Other properties of the set D.}
 
 The next result will be used in \S4, when we  mention the relationship between the set $D$ and the set $E$ (the  definition will be given in below).\\
 
{\bf Proposition 2.}
\begin{theorem}
Let $T_t$ be a linear semigroup in a Banach space $X$, and assume that it has the backward uniqueness property. Put 
$$
	D=\{x\in X|\quad \forall t>0 \quad \exists y=y(t) \quad \hbox{s.t.} \quad 
	T_ty=x \}.
$$
Then, $D$ is a Fr{\'e}chet space, equipped with the following countable norms.
$$
	\|x\|_0=\|x\|,\cdots, \|x\|_n=\|T_{-n}x\|,\cdots \quad n\in {\bf N}\cup \{0\},
$$
where $\|\cdot\|$ denotes the norm in $X$.
\end{theorem}
 $Preoof.$ First, we see that $\|\cdot\|_n$ ($n=0,1,...$) are seminorms. In fact, 
 $\|x\|_n\geq 0$, 
 $$
 	\|\alpha x\|_n = \|T_{-n}(\alpha x)\|=|\alpha|\|x\|_n \quad \forall \alpha\in {\bf R},\quad \forall x\in X,
 $$
 $$
 	\|x+y\|_n = \|T_{-n}(x+y)\|=\|T_{-n}x+T_{-n}y\|
 	\leq \|T_{-n}x\|+\|T_{-n}y\|=\|x\|_n + \|y\|_n.
 $$
Next, we show that $\|\cdot\|_n$ is separable, by a contradiction argument. Assume that there exists $x\neq 0$ such that $\|x\|_n$$=\|T_{-n}x\|$$=0$.  However, this reads $T_{-n}x=0$, and
$x=T_n(T_{-n}x)=0$, which is a contradiction. Thus, $\|x\|_n>0$ for any $x\neq 0$, and $\|\cdot\|_n$ is separable. \\
Finally, we confirm that $D$ is complete with this topology. Let $\{y_m\}\subset D$ be a Cauchy sequence with respect to the above locally compact topology. 
That is, for any $n$, 
$$
	\|y_m-y_{m'}\|_n=\|T_{-n}y_m-T_{-n}y_{m'}\|\to 0 \quad \hbox{as}\quad 
	m,m'\to \infty.
$$
Thus, there exists $z_n\in D$ such that 
$$
	\lim_{m\to \infty}T_{-n}y_m=z_n \quad \hbox{in}\quad X. 
$$
It is easy to see that 
$$
 	z_{n+1}=T_{-1}z_n,
$$ 
for $z_{n-1}$$=\lim_{m\to \infty} T_{-(n-1)}y_m$
$=\lim_{m\to \infty}T_1T_{-n}y_m$$=T_1z_n$. 
This leads to
$$
	T_{n+1}z_{n+1}=T_nT_1z_{n+1}=T_nz_n=z_0 \quad z_0\in D.
$$
We have now $\lim_{m\to \infty}T_{-n}y_m$$=T_{-n}z_0$ in $X$. Therefore, $y_m$ converges to $z_0\in D$, and we have proved the claim.\\

If $D$ is not dense in $X$ (or $\overline{D(A)}$), neither is $D_t$ ($\forall t>0$) dense in $X$ (or $\overline{D(A)}$). \\

{\bf Proposition 3.} 
\begin{theorem}
If there exists $t>0$ such that 
$$
	D_t=\{x\in X|\quad \exists y\in X\quad \hbox{s.t.}\quad T_ty=x\}
$$
is not dense in $\overline{D(A)}$, then 
$$
	\inf\{t|\quad D_t \quad \hbox{is not dense in}\quad \overline{D(A)}\}=0.
$$
\end{theorem}
$Proof.$ We assume that 
$$
	t_0=
	\inf\{t|\quad D_t \quad \hbox{is not dense in}\quad \overline{D(A)}\}\neq 
	0,
$$
and we shall lead to a contradiction. For any $x\in \overline{D(A)}$, and for any $\e>0$, there exists $y\in \overline{D(A)}$ such that $T_{\frac{2t_0}{3}}y\in U_{\e}(x)$. For this $y$, and for any $\delta>0$, there exists $z\in \overline{D(A)}$ such that $T_{\frac{2t_0}{3}}z\in U_{\delta}(y)$. 
Since
$$
	\|T_{\frac{4t_0}{3}}z-x\|\leq \|T_{\frac{4t_0}{3}}z-T_{\frac{2t_0}{3}}y\|
	+\|T_{\frac{2t_0}{3}}y-x\|
$$
$$
	\leq \|T_{\frac{2t_0}{3}}\|\|T_{\frac{2t_0}{3}}z-y\|+\e
	\leq \|T_{\frac{2t_0}{3}}\|\delta+\e,
$$
by tending $\delta\to 0$, we get 
$$
	\|T_{\frac{4t_0}{3}}z-x\|\leq \e.
$$
However, since $x$ is an arbitraly point, this contradicts to the definition of $t_0$. And thus, $t_0=0$ must hold.\\

{\bf Proposition 4.}
\begin{theorem}
	If there exists $t>0$ such that 
$$D_t=\{x\in \overline{D(A)}|\quad \exists y\in \overline{D(A)}\quad 
	\hbox{s.t.}\quad T_ty=x\}
$$
is not dense in $\overline{D(A)}$, then there exists an open set $U$ in $X$ such that 
$$
	U\cap \overline{D(A)}\neq \emptyset,\quad 
	U\cap \overline{D(A)}\cap D_t = \emptyset,
$$
and 
\begin{equation}\label{lr}
	\inf\{t|\quad U\cap D_t=\emptyset\}
	=\min\{t|\quad U\cap D_t=\emptyset \}. 
\end{equation}
\end{theorem}
$Proof.$ Assume that (\ref{lr}) does not hold, and put 
$t_0=$$\inf\{t|\quad U\cap D_t=\emptyset\}$. We have 
$$
	U\cap \overline{D(A)}\cap D_{t_0}\neq \emptyset,
$$
and for any $\e>0$ $U\cap \overline{D(A)}\cap D_{t_0+\e}$ $= \emptyset$.
$T_{-t_0}U$ $=\{x\in \overline{D(A)}|\quad T_{t_0}x\in U\}$ $\neq \emptyset$
 is an open set. Thus, there exists $x\in T_{-t_0}U$ such that $T_{\e}x\in T_{-t_0}U$ for some $\e>0$. However, this reads $T_{t_0+\e}x\in U$, which is a contradicts to the definion of $t_0$. Therefore, (\ref{lr}) must be true.\\

\section{The extended space $E$.}

Assume that there exists a semigroup $T_t$ on a Banach space $X$, and that 
 it defines  solution trajectories in $X$. Assume also that 
 $$
 	D=\{x\in X|\quad \forall t>0 \quad \exists y=y(t)\in X\quad \hbox{s.t.}\quad T_ty=x\}
 $$
is dense in $X$. Then, for an arbitrary point $z\in X$, we can take a sequence of points $\{z_n\}$ in $D$, such that
$$
	\lim_{n\to \infty}z_n=z\quad \hbox{in}\quad X.
$$
Remark that the solution trajectory which passes $z_n$ continues backward until $t=-\infty$. In general, we cannot expect that the solution trajectory which passes $z$ continues until $t=-\infty$. For any $t>0$, we know that $T_{-t}z_n$ exists, but we do not know if the sequence $\{T_{-t}z_n\}$ is a Cauchy sequence or not. Here, our idea is that we regard the sequence $\{T_{-t}z_n\}$ as "a point" in an extended space, and we define 
$$T_{-t}z=\{T_{-t}z_n\}.$$ 
This consideration leads us to define the following extended space $E$, on which the semigroup is extendable to the group. 
\begin{equation}\label{E}
	E=\{
	(z_n)|\quad z_n\in D,\quad \exists t>0 \quad \lim_{n\to \infty}T_tz_n\in X
	/{\sim},
\end{equation}
where $\sim$ represents the equivalence defined by :
$$
	(z_n)\sim (z_n') \quad \hbox{if and only if}\quad \lim_{n\to \infty}T_tz_n
	=\lim_{n\to \infty}T_tz_n' \quad \hbox{for some}\quad t>0.
$$

{\bf Theorem 1.}
\begin{theorem}
	Let $T_t$ be a semigroup defined in a Banach space $X$. Assume that the following two conditions hold.\\

(I) The set 
$$
	D=\{x\in X|\quad \forall t>0 \quad \exists y=y(t)\in X\quad \hbox{s.t.}\quad T_ty=x\}
$$
is dense in $X$.\\
(II) $T_t$ has the backward uniqueness property, i.e. if $T_tx=T_ty$ for some 
$t>0$, then $x=y$.\\

Let $E$ be the set defined in (\ref{E}). Then, the following holds. 
\begin{itemize}
\item $X\subset E$.
\item There exists a group ${\mathcal T_t}$ ($t\in {\bf R}$) on $E$, such that 
$$
	{\mathcal T_t}{\mathcal T_s}={\mathcal T_{t+s}},\quad
	{\mathcal T_0}=I \quad (I\quad \hbox{identity map}).
$$
\item 
$${\mathcal T_t}x=T_tx\quad x\in X\quad (\hbox{whenever}\quad T_tx \quad\hbox{exists}).
$$
\end{itemize}
\end{theorem}
$Proof.$ The relationship $X\subset E$ comes from the identification of $x\in X$ to $(x_n)$$\in E$, where the sequence $\{x_n\}$ in $X$ satisfies $\lim_{n\to \infty}x_n$$=x$ in $X$. We define the group ${\mathcal T_t}$ in the following way : for $(x_n)\in E$
$$
	{\mathcal T_t}(x_n)=(T_tx_n) \quad \forall t\in {\bf R},
$$
where the right-hand side is well-definded from the boundedness of $T_t$. 
Then, 
$$
	{\mathcal T_{t+s}}(x_n)=(T_{t+s}x_n) =(T_tT_sx_n)
	={\mathcal T_t}{\mathcal T_s}(x_n);
$$
$$
	{\mathcal T_0}(x_n)=(x_n),\quad \hbox{if}\quad 
	\lim_{n\to \infty}x_n=x, \quad\hbox{then}\quad 
	{\mathcal T_t}(x_n)=(T_t x). 
$$
\medskip

The semigroup $T_t$ on $X$ is extended to the group ${\mathcal T_t}$ on $E$, 
and the solution trajectory defined by $T_t$ is extended to ${\mathcal T_t}$ on $E$, backward in time until $t=-\infty$. \\

{\bf Remarks} 1. For the case that $T_t$ is a nonlinear semigroup, the above theorem holds by replacing $X$ to $D(A)$, too.\\
2. Even if $T_t$ does not have the backward uniqueness property, ${\mathcal T_t}$ constructed as above has the backward uniqueness property on $E$ defined by (\ref{E}). This is the reason why we assumed the condition (ii).\\
3. If a semigroup $T_t$ corresponds to a group ${\mathcal T_t}$ by ${\mathcal T_t(x_n)}=(T_tx_n)$, the infinitesmal generator $A$ of $T_t$ corresponds to the infinitesmal generator ${\mathcal A}$ defined by 
$$
	{\mathcal A}(x_n)=(Ax_n)\quad \hbox{whenever the right-hand side exists.}
$$

\section{The relationship between  $E$ and $D$.}

In this section, we study some structures of $E$ introduced in \S3, and the relationship between $E$ and $D$, the set of initial points from which the solution trajectory continues until $t=-\infty$ in $X$. \\

The space $E$ is not in general a Banach space, but the intermediate spaces between $X$ and $E$ are Banach spaces.\\

{\bf Proposition 5.} 
\begin{theorem} Let $T_t$ ($t\geq 0$) be a semigroup on a Banach space. 
	For $t>0$, put 
$$
	E_{-t}=\{(z_n)|\quad z_n\in D\quad \lim_{n\to \infty} T_t z_n\in X\},
$$
the set of points in $X$ from which the solution trajectory exists backward in time at least for $t$. Then, $E_{-t}$ is a Banach space equipped with the norm : $$
	\|(z_n)\|_{-t}=\|\lim_{n\to \infty}T_tz_n\|. 
$$
\end{theorem}
$Proof.$ It is clear that $E_{-t}$ is a linear space. To see that $E_{-t}$ is a normed space, we confirm the following. \\
$$
	\|(z_n)+(y_n)\|_{-t}=\|\lim_{n\to \infty}T_t(z_n+y_n)\|
		\qquad\qquad\qquad\qquad
$$
$$\qquad\qquad
	\leq \|\lim_{n\to \infty}T_tz_n\|+ \|\lim_{n\to \infty}T_ty_n\|
	= \|(z_n)\|_{-t}+\|(y_n)\|_{-t}.
$$
Also, it is obvious that $\|(z_n)\|_{-t}\geq 0$, and that if $\|(z_n)\|_{-t}=0$, then $(z_n)=0$, for $\|\lim_{n\to \infty}T_tz_n\|=0$ implies $\lim_{n\to \infty}T_tz_n=0$ and thus $(z_n)\sim(0)$. \\
Moreover, for any $\alpha \in {\bf R}$, 
$$
	\|\alpha (z_n)\|_{-t}=\|(\alpha z_n)\|_{-t}=\|\lim_{n\to \infty}T_t(\alpha z_n)\|
		=|\alpha |\|(z_n)\|_{-t}.
$$
Finally, to see that $E_{-t}$ is complete with respect to the norm $\|\cdot\|_{-t}$, let $\{(z_{n}^{m})_m\}_m$ be a Cauchy sequence with respect to the norm $\|\cdot\|_{-t}$, i.e. 
$$
	\|(z_{n}^m)_m-(z_{n}^{m'})_{m'}\|_{-t}\to 0\quad \hbox{as}\quad m,m'\to +\infty.
$$
If we put $c_m=\lim_{n\to \infty}T_tz_{n}^m=c_m$, then
$$
	\|c_m-c_m'\|\to 0\quad \hbox{as}\quad m,m'\to +\infty.
$$
Since $\{c_m\}$ is a Cauchy sequence with respect to the norm $\|\cdot\|$, there exists $c_0\in X$ such that $\lim_{m\to \infty}c_m$$=c_0$ in $X$. Since $D$ is dense in $X$, there exists $\{y_n\}$$\subset D$ such that $\lim_{n\to \infty}y_n=c_0$ in $X$. Therefore, 
$$
	\|(z_{n}^m)-(T_{-t}y_n)\|_{-t}=\|c_m-c_0\|\to 0\quad m\to \infty.
$$
From the above argument, we have shown that $E_{-t}$ is a Banach space.\\

{\bf Remarks.} 1. Similarly, we can see that the intermediate spaces between D and E :
$$
	D_t=T_tX=\{x\in X|\quad \exists y\in X,\quad T_ty=x\}
$$
is also a Banach space with respect to the norm $\|x\|_t=\|T_{-t}x\|$. \\
2. We have the following relationship :
$$
	D \subset D_t \subset X \subset E_{-t} \subset E.
$$
3. Let $S_t$ be the semigroup defined in Corollary 1 (\S2.1), i.e.
$$
	S_t x_0=x(t)=T_tx_0+\int_0^t T_{t-s}f(s)ds,
$$
where $T_t=e^{-tA}$. Define the norm of $E_{-t}$ by
$$
	\|x\|_{-t}=\|S_tx-\int_0^t T_{t-s}f(s)ds\|.
$$
Then, $E_{-t}$ is a Banach space with respect to the norm $\|\cdot\|_{-t}$. \\
4. Let $T_t$ be a non-expansive nonlinear semigroup defined on $\overline{D(A)}$, a closed convex subset in a Banach space $X$. Assume that the set 
$$
	D=\{x\in \overline{D(A)}|\quad \forall t>0 \quad \exists y=y(t)\quad 
	\hbox{s.t.}\quad T_ty=x\}
$$
is dense in $\overline{D(A)}$. Put 
$$
	\overline{D(A)}_{-t}=
	\{(z_n)|\quad \lim_{n\to \infty}T_tz_n\in \overline{D(A)},\quad z_n\in D\}.
$$
Then, $\overline{D(A)}_{-t}$ is a complete metric set with respect to the norm : 
$$
	\|(z_n)\|_{-t}=\|\lim_{n\to \infty}T_t z_n\|. 
$$
However, in general, we do not know whether it is possible to linearize $\overline{D(A)}_{-t}$, and extend it to a larger Banach space containing the original space $X$. It is because, differently from the linear case, there is no longer the structure of a Banach space. \\

\medskip

Next, by assuming that $X=H$ a Hilbert space, we consider the relationship between $E$ and $D$. \\

{\bf Definition 1.} We say that two linear vector spaces $V$ and $W$ are dual, with respect to $<,>_{V\times W}$, when the following holds.
$$
	<x,y>_{V\times W}=0\quad \forall y\in W\quad\hbox{implies}\quad x=0.
$$

{\bf Proposition 6.}
\begin{theorem}
The linear vector spaces $V$ and $W$ are dual, with respect to the following 
 bilinear functional : for any $x\in D$ and any $(z_n)\in E$,
 $$
 	<x,(z_n)>_{D\times E}=(T_{-t}x,\lim_{n\to \infty}T_t z_n),
 $$
 where $(\cdot,\cdot)$ is the inner product of $H$, and 
 $$
 	t=\inf\{s|\quad \lim_{n\to \infty}T_sz_n\in H\}. 
 $$
\end{theorem}
$Proof.$ For any $(z_n)\in E$, let $t=\inf\{s|\quad $$\lim_{n\to \infty}T_sz_n$$\in H\}$, and assume that 
$$
	(T_{-t}x,\lim_{n\to \infty}T_t z_n)=0.
$$
In particular, if for any $(z)\in E$ such that $z\in D$ the above holds, then from the density of $D$ in $H$, we get $x=0$. Thus, we proved the claim.\\

{\bf Remark.} In Proposition 2, we have seen that $D$ is a Fr{\'e}chet space with respect to the countable seminorms $\|\cdot\|_n$ ($\|x\|_n$$=\|T_{-n}x\|$). Concerning with this topology, the  linear functional :
$$
	<\cdot,(z_n)>_{D\times E}=(T_{-t}\cdot,\lim_{n\to \infty}T_tz_n)
$$
is an element in $D^{\ast}$. This shows the existence of a one to one map 
from $E$ to $D^{\ast}$. 

Finally, let us introduce a topology by the countable number of norms : $\|x\|_n$$=\|T_{-n}x\|$. \\

{\bf Proposition 7.}
\begin{theorem}
	Let $H$ be a Hilbert space, and $\{{\bf e_n}\}$ be a complete orthonormal bases of $H$. If a semigroyp $T_t$  on $H$ has the expression 
$$
	T_t u_0=\sum_{n=1}^{\infty} \exp(\lambda_n t) a_n {\bf e_n} \quad t\geq 0,
	\quad u_0=\sum_{n=1}^{\infty} a_n{\bf e_n}, 
$$
then $D^{\ast}=E$. 
\end{theorem}
$Proof.$ It is clear that $E\subset$$D^{\ast}$. For $f\in D^{\ast}$, put $b_k=<{\bf e_k},f>_{D\times D^{\ast}}$, and see 
$$
	\sum_{n=1}^{\infty}a_nb_n<\infty \quad \hbox{if}\quad \sum_{n=1}^{\infty}a_n{\bf e_n}\in D.
$$
Since the sequence $(a_n)$ satisfies $(a_n\exp\l_n t)$$\in (l^2)$ for any $t\in {\bf R}$, we can write 
$$
	b_n=\beta_n\exp \lambda_n t\quad (\beta_n)\in (l^2),
$$
where 
$$
	t=\sup\{s|\quad\exists (\beta_n)=(\beta_n(s))\in (l^2)\quad
    \hbox{s.t.}\quad b_n=\beta_n\exp\lambda_n s\}.
$$ 
Here, $(\sum_{k=1}^{n}\beta_k \exp \lambda_k t{\bf e_k})\in E$. If $t\geq 0$, 
$$
	<{\bf e_l},(\sum_{k=1}^{n} \beta_k\exp \lambda_k t {\bf e_k})>_{D\times E}
		=\beta_l \exp \lambda_l t=b_l,
$$
 and  if $t< 0$, 
$$
	<{\bf e_l},(\sum_{k=1}^{n} \beta_k\exp \lambda_k t {\bf e_k})>_{D\times E}
		=<\exp \lambda_l t {\bf e}_l, \sum_{k=1}^{\infty}\beta_k {\bf e_k}>=\beta_l \exp \lambda_l t=b_l.
$$
    Thus, 
$$
    f=(\sum_{k=1}^{n}\beta_k \exp \lambda_k t {\bf e_k})\in E,
$$
and we proved the claim. \\


\begin{thebibliography}{31}
\bibitem{br1}
H. Brezis. Operateurs Maximaux Monotones et Semi-Groupes de Contractions Dans Les Espaces De Hilbert. North Holland (1973).
\bibitem{br2}
H. Brezis. Analyse fonctionnelle Theorie et applications. Masson (1987).
\bibitem{henry}
D. Henry. Geometric Theory of Semilinear Parabolic Equations. Springer. Lecture Notes in Mathematics. \textbf{840} (1981).
\bibitem{mp}
M. Pierre.  Enveloppe d'une famille de semi-groupes non lineaires et equations d'evolution. Seminaire d'analyse non lineaire. Universite de Besancon. (1976-1977). 
\end{thebibliography}
\end{document}